\input AHTOHFIE.STY
\headline{\ifnum\count0=1 \vtt UDC 512.543.7+512.543.16\hss\else\hss\fi}
\centerline{\ssdbf
THE SQ-UNIVERSALITY OF ONE-RELATOR RELATIVE PRESENTATIONS%
\footnote{}{\rm
This work was supported by the Russian Foundation for Basic Research,
project no. 05-01-00895.} }

\smallskip
\centerline{\ss Anton A. Klyachko}
\smallskip
{
\ssqi
\centerline{Faculty of Mechanics and Mathematics, Moscow State University}
\centerline{Moscow 119992, Leninskie gory, MSU}
\centerline{klyachko@daniil.math.msu.su}
}

\medskip

\abstract{%
Adding two generators and one arbitrary relator to a nontrivial
torsion-free group, we always obtain an SQ-universal group. In the course
of the proof of this theorem, we obtain some other results of independent
interest. For instance, adding one generator and one relator in which the
exponent sum of the additional generator is one to a free product of two
nontrivial torsion-free groups, we also obtain an SQ-universal group.
\hfil\break
{\itsmall Key words}:
relative presentations, one-relator groups, SQ-universality, equations
over groups.
\hfil\break
{\itsmall MSC}\/:
20E06, 20F05, 20F06.
}

\s 1.
Introduction

Recall that a group $G$ is called {\it SQ-universal} if every
countable group can be embedded in some quotient of~$G$.
Examples of SQ-universal groups are all nonabelian
free groups; all free products, except the infinite dihedral group
$\Z_2*\Z_2$ (see [LS77]); many amalgamated free products and
HNN-extensions ([Lo86], [LS77]); all finite-index subgroups of
SQ-universal groups and all virtually SQ-universal groups [Neu73]; all
nonelementary hyperbolic groups [Ols95] and even all (with some obvious
exceptions) relatively hyperbolic groups, in particular, all groups with
infinitely many ends~[AMO06].

The starting point for our study is the following theorem.

\proclaim{Sacerdote--Schupp theorem \rm [SaSc74] (see also [LS77])}.
A group admitting a presentation with one relator and at least three
generators is SQ-universal.

Actually, there is a much more general fact.

\proclaim{Baumslag--Pride theorem \rm [BaPr78]}.
A group having a presentation with two more generators than relators
is SQ-universal.
Moreover, any such group is {\it large in the sense of Gromov}, i.e. it
has a finite-index subgroup admitting an epimorphism onto a
nonabelian free group.

A further generalisation of the Sacerdote--Schupp theorem is the following
result.

\proclaim{St\"ohr--Gromov theorem \rm[St83], [Gr83]}.
A group having a presentation in which there are more generators
than relators and one of the relators is a
proper power is SQ-universal and even large in the sense of Gromov.%
\fn{%
A {\it proper power} is an element of a free group $F$ of the form
$u^k$, where $u\in F$ and $\Z\ni k\ge2$. In particular, the identity
element is a proper power. Thus,
the St\"ohr--Gromov theorem is a generalisation
of the Baumslag--Pride theorem.
}

Further results on this subject can be found in, e.g.,
[Ed84],
[How98],
[Bu05],
[La05],
[OlOs06].
In this paper, we generalise the Sacerdote--Schupp theorem in another
direction.

Let $G$ be a group. A group given by a one-relator relative presentation
over $G$ is
$$
\~G=\pres<G,x_1,x_2,\dots,x_n|w=1>\:=
G*F(x_1,x_2,\dots,x_n)/\nc w.
$$
Here $x_1,\dots,x_n$ are some letters (not belonging to $G$)
and $w$ is a word in the alphabet
$G\cup\{x_1^{\pm1},\dots,x_n^{\pm1}\}$ (such a word can be considered as
an element of the free product $G*F(x_1,x_2,\dots,x_n)$ of
$G$ and the free group with basis $x_1,x_2,\dots,x_n$). In other words,
the presentation of the group $\~G$ is obtained from a presentation
$G=\pres<A|R>$ of $G$ by adding several new generators and one
new relator:
$\~G=\pres<A\cup\{x_1,x_2,\dots,x_n\}|R\cup\{w\}>$.

\Th 1.
If $G$ is a nontrivial torsion-free group and $n\ge2$, then the group
$
\~G=\pres<G,x_1,x_2,\dots,x_n | w=1>
$
is SQ-universal for any $w\in G*F(x_1,\dots,x_n)$.

\Corollary {\rm [Kl06b]}.
Under the conditions of Theorem 1 the group $\~G$ {\rm (as well as any
SQ-universal group)} has a nonabelian free subgroup.

\Remark 1.
It is easy to show that the group $\~G$ from Theorem 1 need not be
large in the sense of Gromov.

\Remark 2.
Certainly, the assertion of Theorem 1 is not valid for $n=1$
(as well as the assertion of the Sacerdote--Schupp theorem is not
valid for groups with two generators). However, the group $\~G$
with $n=1$ has some properties weaker than the SQ-universality.
In particular,  $\~G$ is nontrivial [Kl93], it cannot be a nonabelian
simple group (if $w\ne g_1x_1^{\pm1}g_2$) [Kl05], and the natural
mapping $G\to\~G$ is nonsurjective (if $w\ne g_1x_1^{\pm1}g_2$) [CR01];
in the case when the exponent sum of $x_1$ in $w$ is
$\pm1$, the group $\~G$ always (with some obvious exceptions)
contains a nonabelian free subgroup [Kl06b].

Theorem 1 considers relative presentations with at least
two additional generators; however, an important role in the proof
of this theorem is played by the study of one-generator relative
presentations
$$
\~G=\gp{G,t\ |\ w=1}\:=(G*\gp t_\infty)/\nc w,
\quad\hbox{where }
w\equiv \prod g_it^{\epsilon_i},
\ g_i\in G,\ \epsilon_i\in\{\pm1\}.
\eqno{(1)}
$$
Such a presentation is called {\it unimodular} if
$\sum\epsilon_i=\pm1$. It is known that unimodular relative
presentations have some good properties and are more convenient to
study (see, e.g., [Kl93], [Kl94], [FeR96], [CR01], [FoR05],
[Kl05], [Kl06a], [Kl06b]).

In [Kl06a], we suggested a generalisation of the notion of
unimodularity to the so-called
{\it generalised relative presentation}
$$
\~G=\pres<G*T|\prod_{i=1}^n g_it_i=1>\:=(G*T)\Big/
\!\!\gp{\!\!\gp{\prod g_it_i}\!\!}.
\eqno{(*)}
$$
Here $T$ is a group (not necessarily cyclic),
$g_i\in G$, and $t_i\in T$; the word $\prod\limits_{i=1}^n g_it_i$ is
assumed to be cyclically reduced.

The generalised relative presentation $(*)$ over a group $G$
is called {\it unimodular} if
\item{1)} the order of the element $\prod t_i$ is infinite
in the group $T$;
\item{2)} the cyclic subgroup $\gp{\prod t_i}$ is normal in $T$;
\item{3)} the quotient group $T/\gp{\prod t_i}$ possesses the
          strong unique-product property.

\noindent
Recall that a group $H$ is called a {\it UP-group, {\rm or a group with
the} unique product property}, if the product $XY$ of any two finite
nonempty subsets $X,Y\subseteq H$ contains at least one element which
decomposes uniquely into the product of an element from $X$ and an element
from $Y$.%
\fn{Some time ago, there was the conjecture that any torsion-free
group is UP (the converse is, obviously, true).  However, it turned out
that there exist counterexamples ([P88], [RS87]).}

We say that a group $H$ has the {\it strong unique product property} if
the product $XY$ of any two finite
nonempty subsets $X,Y\subseteq H$ such that $|Y|\ge2$ contains at least
two uniquely decomposable elements $x_1y_1$ and $x_2y_2$ such that
$x_1,x_2\in X$,\ \ $y_1,y_2\in Y$, and $y_1\ne y_2$.

As far as we know, all known examples of UP-groups have the strong
UP-property. In particular, all right orderable groups, locally indicable
groups, and diffuse groups in the sense of Bowditch have the strong UP
property.

Theorem 1 is an easy corollary of the following theorem.

\Th 2.
If the generalised relative presentation $(*)$ over a noncyclic
torsion-free group $G$ is unimodular and the group $T$ is not cyclic, then
the group $\~G$ given by presentation $(*)$ is SQ-universal.

In order to prove Theorem 2, we establish the following fact about
usual (nongeneralised) unimodular relative presentations.

\Th 3.
If $G_1,\dots,G_l$ are noncyclic torsion-free groups,
$l\ge2$, and a relative presentation
$$
L=\gp{G_1*\dots*G_l,t\ |\ w=1}
$$
over
the group $G_1*\dots*G_l$ is unimodular, then the group $L$ is SQ-universal.
Moreover, each countable group $S$ embeds into a quotient group
$L/N$, in which the Freiheitssatz holds, i.e.,
$$
\gp{t,G_{i_1},\dots,G_{i_{l-1}}}=\gp t_\infty*G_{i_1}*\dots*G_{i_{l-1}}
\quad\hbox{in }L/N
$$
if the word $w$ is conjugate in group $\gp t_\infty*G_1*\dots*G_l$
to no element of the subgroup $\gp t_\infty*G_{i_1}*\dots*G_{i_{l-1}}$.

According to [Kl06a], we say that presentation (1) is
{\it magnusian}
if the natural mapping $G\to\~G$ is injective and
$\gp{H,t}=H*\gp{t}_\infty$ in the group $\~G$
(i.e., $t$ is {\it transcendental} over
$H$ in $\~G$) for any free factor $H$ of $G$ such that
$w$ is not conjugate in $G*\gp{t}_\infty$ to an element of
$H*\gp{t}$.

In [Kl06a], we showed that every unimodular presentation
over a torsion-free group is magnusian. To prove the main
results of this paper, we need a stronger property of
unimodular presentations.

We say that presentation (1) is
{\it strongly magnusian} if
the element $t$ is transcendental in $\~G$ over each subgroup
$H\subseteq G$ such that
\item{1)}
$w$ is not conjugate in $G*\gp{t}_\infty$ to an element of $H*\gp{t}$;
\item{2)}
each coefficient $g_i$ either lies in $H$ or is transcendental over $H$.

\Proposition 1.
If presentation (1) is unimodular and each nonidentity coefficient $g_i$
has infinite order in $G$, then presentation (1) is strongly
magnusian.

The proof of Proposition 1 is based on several earlier known
results, in particular, on the following theorem.

\Th 4.
Suppose that $A\zvezda_C B$ is the free amalgamated product of groups $A$
and $B$ with amalgamated subgroup~$C$, $v=b_0a_0\dots b_ma_mb_{m+1}\in
(A\zvezda_C B)$, $m\ge1$, and each coefficient of the word $v$ (except
maybe the first and the last ones) is transcendental over $C$, i.e.,
$\gp{a_i,C}=\gp{a_i}_\infty*C$ in $A$
for $i=0,\dots,m$
and
$\gp{b_i,C}=\gp{b_i}_\infty*C$ in $B$
for $i=1,\dots,m$.
Then, for any automorphism~$\phi$ of the group $B$, the natural mappings
$$
A\to\gp{A\zvezda_C B\ \biggm|\ \{b^v=b^\phi\ |\ b\in B\}} \leftarrow B
$$
are injective.

This theorem was proven in [Kl94], but it has never been published. The
last section of this paper contains a proof of Theorem 4. Contrary to
purely algebraic arguments of all preceding sections, the proof of Theorem
4 is geometric. Other known facts from which Proposition 1 is derived
have similar proofs.

\proclaim{Notation}\rm{
which we use is mainly standard.
Note only that if $k\in\Z$, $x$ and $y$ are elements of a
group, and $\phi$ is a homomorphism from this group into another
group, then $x^y$, $x^{ky}$, $x^{-y}$, $x^\phi$, $x^{k\phi}$, and
$x^{-\phi}$ denote $y^{-1}xy$, $y^{-1}x^ky$, $y^{-1}x^{-1}y$, $\phi(x)$,
$\phi(x^k)$, and $\phi(x^{-1})$, respectively;
the commutator $[x,y]$ is understood as $x^{-1}y^{-1}xy$.
If $X$ is a subset of a group, then
$\gp{X}$ and $\nc{X}$ denote the subgroup generated by $X$
and the normal subgroup generated by the set $X$, respectively.
The symbol $|X|$ denotes the cardinality of the set $X$}.


\s 2.
Proof of Theorem 1

If the group $G$ is cyclic, then $\~G$ is a group with
at least
three generators and one relator. The SQ-universality of such groups
is asserted by the Sacerdote--Schupp theorem. Thus, we assume that
$G$ is noncyclic.

Suppose that the word $w$ has the form
$w\equiv g_1x_{j_1}^{\epsilon_1}g_2x_{j_2}^{\epsilon_2}\dots
g_px_{j_p}^{\epsilon_p}$ and the word $w'\in F(x_1,\dots, x_n)$ is
obtained from~$w$ by erasing the coefficients:
$w'=x_{j_1}^{\epsilon_1}x_{j_2}^{\epsilon_2}\dots x_{j_p}^{\epsilon_p}$.

\smallskip
\noindent{\bf Case 1}:
$w'$ is a proper power in the free group $F(x_1,\dots, x_n)$.  In this
case, the group $\~G$ is SQ-universal, because the one-relator homomorphic
image $T_1=\gp{x_1,\dots,x_n\ |\ w'=1}$ of $\~G$ is SQ-universal by the
St\"ohr--Gromov theorem.

\smallskip
\noindent{\bf Case 2}:
$w'$ is not a proper power.
Consider the groups
$$
T=\gp{x_1,\dots,x_n\ |\ [x_1,w']=\dots=[x_n,w']=1}
\quad\hbox{and}\quad
T_1=\gp{x_1,\dots,x_n\ |\ w'=1}=T/\gp{w'}.
$$
The group $T$ is the free central extension of the one-relator group
$T_1$. It is well known that if $w'$ is not a proper power
in the free group $F(x_1,\dots,x_n)$, then the group $T_1$ is locally
indicable ([B84]) and, therefore, has strong unique-product property.
The element $w'$ has infinite order in $T$ (see [LS77]). Thus,
the generalised relative presentation $\pres<G,T|w=1>$ is unimodular.
The group $T$ is not cyclic, because its commutator quotient
is the free abelian group of rank $n\ge2$.  Therefore, by Theorem 2, the
group $\pres<G,T|w=1>$ is SQ-universal. It remains to note that this group
is a homomorphic image of $\~G$, and a group having SQ-universal
homomorphic image is SQ-universal itself.


\s 3.
Iterated amalgamated free products

In this section, we reproduce a construction from [Kl06a]
in a more general situation.

Let $I$ be a set, and let $\Omega$ be a
family of subsets of $I$. For each $i\in I$, let $G_i$
be a group, and for each
$\omega\in\Omega$, let $G_\omega$ be a quotient of the free product
$\zvezda\limits_{i\in\omega}G_i$:
$$
G_\omega= \left(\zvezda_{i\in\omega}G_i\right)\Big/N_\omega.
$$
The natural question arises: under what conditions are the natural mappings
$$
\phi_\omega\:G_\omega\to G_I\:=
\left(\zvezda_{i\in I}
G_i\right)\Bigg/\!\gp{\!\!\!\gp{\bigcup_{\omega\in\Omega}N_\omega}\!\!\!}
$$
injective? Or under what conditions can the group $G_I$ be
considered as an amalgamated free product of the groups~$G_\omega$?

The following proposition gives some sufficient condition for this
question to have a positive answer.

\Proposition 2.
Suppose that
$$
N_\omega\cap\zvezda\limits_{j\in\omega\setminus\{i\}}G_j=\1
\eqno{(**)}
$$
for each $\omega\in\Omega$ and each
$i\in\omega\setminus\left(\bigcap\Omega\right)$. Suppose also that,
for each finite subfamily $F\subseteq\Omega$ with
$|F|\ge2$, there exist elements $\min, \max\in \bigcup F$ such
that
\item{\rm 1)}
the element $\min$ belongs to precisely one set
$\omega_{\min}\in F$;
\item{\rm 2)}
the element $\max$ belongs to precisely one set
$\omega_{\max}\in F$;
\item{\rm 3)}
$\omega_{\min}\ne\omega_{\max}$.
\enditem
Then all of the natural mappings $\phi_\omega\:G_\omega\to G_I$ are
injective.

\Example.
Suppose that $I=\{a,b,c,d,e,f\}$ and
$\Omega=
\left\{
\{a,b,d,e\},
\{b,c,e,f\},
\{d,e,f\}
\right\}$.

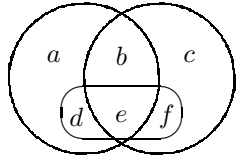
\nobreak
\centerline{\input 0.PIC}
\goodbreak

\noindent
Let $A,\dots,F$ be
the corresponding six groups $G_i$, and let
$\bf ABDE$, $\bf BCEF$, and $\bf DEF$ be the three groups
$G_\omega$.
It is easy to see that conditions 1), 2), and 3) hold for
the family $\Omega$ and each of its two-set subfamilies. Suppose that
condition~$(**)$ holds too. Then the validity of Proposition 2 (for
this example) is implied by the following decomposition of $G_I$ into
an amalgamated free product:
$$
G_I= \left( ({\bf DEF}*B)
\zvezda_{B*D*E} {\bf ABDE} \right) \zvezda_{B*E*F} {\bf BCEF}.
$$

\smallskip

To prove Proposition 2 in the general case, we need a lemma.

\Lemma 1.
Suppose that the conditions of Proposition 2 hold,
$\Omega'$ is a finite subfamily of $\Omega$,
$\omega\in\Omega$, and
$\alpha\subseteq\omega\cap\left(\bigcup\Omega'\right)$
is a proper subset of $\omega$ contained in
$\bigcup\Omega'$ and containing $\bigcap\Omega$. Then the natural mapping
$$
\zvezda_{i\in\alpha}G_i\to G_{\Omega'}\:=
\left(\zvezda_{i\in\bigcup\Omega'}
G_i\right)\Bigg/
\!\gp{\!\!\!\gp{\bigcup_{\omega'\in\Omega'}N_{\omega'}}\!\!\!}
$$
is injective.

\Proof

\noindent{\bf Case 1: $\omega\in\Omega'$.}
Let us use induction on the cardinality of $\Omega'$.
If $|\Omega'|=1$ (i.e., $\Omega'=\{\omega\}$), then the assertion of
Lemma~1 is true by condition $(**)$. Suppose that $|\Omega'|\ge2$. In this
case, according to conditions 1), 2), and 3), the family $F=\Omega'$
contains a set $\omega'\ne\omega$ that contains an element $m\in\omega'$
not lying in $\bigcup(\Omega'\setminus\{\omega'\})$.

By the induction hypothesis (applied to the set $\omega'$ as $\omega$
and the family $\Omega'\setminus\{\omega'\}$ as~$\Omega'$),
the groups
$$
G_i\ \hbox{with}\
i\in\beta\:=\omega'\cap\left(\bigcup(\Omega'\setminus\{\omega'\})\right)
$$
freely generate their free product in the group
$G_{\Omega'\setminus\{\omega'\}}$.  But according to condition $(**)$, the
same groups $G_i$ with $i\in\beta$ freely generate their free product in
the group $G_{\omega'}$ (because $\omega'$ contains an element $m$ not
lying in $\beta$).  Therefore, the group $G_{\Omega'}$ decomposes into the
amalgamated free product of $G_{\Omega'\setminus\{\omega'\}}$ and
$G_{\omega'}$ with amalgamated subgroup $\zvezda_{i\in\beta}G_i$. The
groups $G_i$ with $i\in\alpha$ lie in the factor
$G_{\Omega'\setminus\{\omega'\}}$. Therefore, the assertion of Lemma 1
follows from the induction hypothesis applied to the set $\omega$ and the
family $\Omega'\setminus\{\omega'\}$ as $\Omega'$.

\noindent{\bf Case 2: $\omega\notin\Omega'$.}
In this case, the proof is similar.
We again use induction on the cardinality of $\Omega'$.
If $\Omega'=\emptyset$, then we have nothing to prove.
Suppose that $|\Omega'|\ge1$. In this
case, according to conditions 1), 2), and 3), the family
$F=\Omega'\cup\{\omega\}$
contains a set $\omega'\ne\omega$ with an element $m\in\omega'$
not lying in $\bigcup(F\setminus\{\omega'\})$ (see Fig. 1).

\goodbreak
\bigskip
\centerline{\input 1.PIC}
\nobreak%
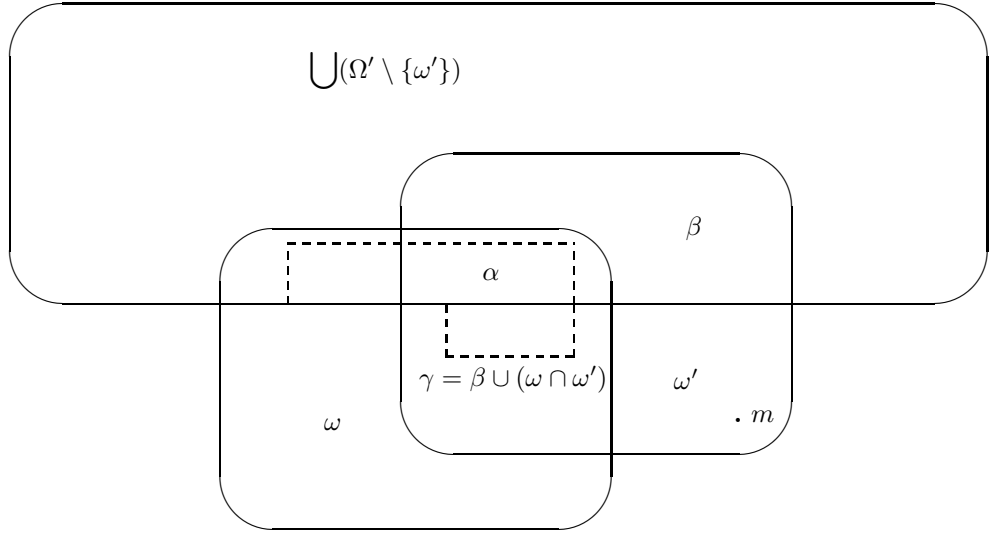
\centerline{Fig. 1}%
\goodbreak
\bigskip

By the induction hypothesis (applied to the set $\omega'$ as $\omega$ and
the family $\Omega'\setminus\{\omega'\}$ as $\Omega'$), the groups
$$
G_i\ \hbox{with}\
i\in\beta\:=\omega'\cap\left(\bigcup(\Omega'\setminus\{\omega'\})\right)
$$
freely generate their free product in the group
$G_{\Omega'\setminus\{\omega'\}}$.
Therefore, the groups
$$
G_i\ \hbox{with}\
i\in\gamma\:=\beta\cup(\omega\cap\omega')=
\omega'\cap\left(\bigcup((\Omega'\cup\omega)\setminus\{\omega'\})\right)
$$
freely generate their product in the group
$$
H=\left(\zvezda_{j\in(\omega\cap\omega')\setminus\beta}G_j\right)*
G_{\Omega'\setminus\{\omega'\}}.
$$
But condition $(**)$ implies that the same groups $G_i$ with $i\in\gamma$,
freely generate their free product in $G_{\omega'}$ (because $\omega'$
contains an element $m$ not lying in $\gamma$). Therefore, the group
$G_{\Omega'}$ decomposes into the amalgamated free product of the groups
$H$ and $G_{\omega'}$:
$$
G_{\Omega'}=H \zvezda_{\gp{G_i\ ;\ i\in\gamma}} G_{\omega'}.
$$
The groups $G_i$ with $i\in\alpha$ lie in the
factor $H$. Therefore, by the induction hypothesis applied to the set
$\omega$ and the family $\Omega'\setminus\{\omega'\}$ as $\Omega'$, the
groups $G_i$ with
$i\in\alpha\cap\left(\bigcup(\Omega'\setminus\{\omega'\}\right)$
freely generate their free product in
$G_{\Omega'\setminus\{\omega'\}}$. This immediately implies that
the groups
$G_i$ with subscripts $i\in\alpha$ freely generate their free product
in $H$ and, hence, in the group $G_{\Omega'}$ containing
$H$ as a subgroup. Lemma 1 is proven.

\smallskip\noindent
{\bf Proof of Proposition 2.}
Clearly, it is sufficient to prove Proposition 2 for a finite family
$\Omega$ of cardinality larger than one. In this case,
$$
G_I=G_\Omega*\left(\zvezda_{i\notin\bigcup\Omega}G_i\right),
$$
and $G_\Omega$ decomposes into the amalgamated free product:
$$
G_\Omega=G_{\omega_{\min}}\zvezda_K G_{\Omega\setminus\{\omega_{\min}\}},
$$
where the amalgamated subgroup $K$ is (by virtue of Lemma 1) the free
product of the groups $G_i$ with
$i\in\omega_{\min}\cap\bigcup(\Omega\setminus\{\omega_{\min}\})$.
An obvious inductive argument completes the proof.


\s 4.
Proof of Theorem 2. The nonsplitting case

In this section we prove Theorem 2 in the case when the
group $\gp{\{t_i\}}\subseteq T$ is noncyclic.

Take an arbitrary countable group $S$.
Put $t=\prod t_i$ and let us decompose $T$ into the union of cosets:
$$
T=\coprod_{x\in T/\gp{t}} c_x\gp{t}, \quad\hbox{where }c_1=1.
$$
For each $x\in T/\gp{t}$, consider an isomorphic copy $G^{(c_x)}$ of the
group $G$ assuming that the isomorphism maps $g\in G$ to
$g^{(c_x)}\in G^{(c_x)}$.
Let us rewrite the relation $\prod g_it_i=1$ in the form
$$
t\prod_i {g_i}^{c_{x_i}t^{k_i}}=1.
\eqno{(2)}
$$
Let $X_1$ be the set of all $x\in T/\gp{t}$
occurring in the reduced form of relation (2). Note that
$|X_1|\ge2$, because $\gp{\{t_i\}}\ne\gp t$ in the case
under consideration.
Put
$$
H_1=\zvezda_{y\in X_1}G^{(c_y)}
$$
and consider the unimodular relative presentation
$$
\~H_1=\pres<H_1,z|z\prod_i {g_i}^{(c_{x_i})z^{k_i}}=1>
$$
over the group $H_1$. Theorem 3 implies
that $\~H_1$ has a quotient group
$K_1$ such that
\item{1)} the group $S$ embeds into $K_1$;
\item{2)} in the group $K_1$, we have the decomposition
$$
\gp{z,\{G^{(c_y)}\ ;\ y\in Y\}}=
\gp z_\infty*\left(\zvezda_{y\in Y} G_y\right)
\eqno{(3)}
$$
for each proper subset $Y\subset X_1$.
\enditem
The group $K_1$ is a quotient of the group
$$
L_1=H_1*\gp z_\infty=
\left(\zvezda_{y\in X_1}G^{(c_y)}\right)*\gp z_\infty
$$
by a normal subgroup $N_1$.

Now, consider the free product
$$
L=\left(\zvezda_{y\in T/\gp t}G^{(c_y)}\right)*\gp z_\infty.
$$
The group $T$ acts on the right on the group $L$ by automorphisms:
$$
z^x           = z^{\epsilon_x},\
\left(g^{(c_y)}\right)^x = g^{(c_{yx})z^l},
$$
where $x\in T$, $y\in T/\gp t$, $\epsilon_x=\pm1$ depending on whether or
not $x$ and $t$ commute, and the integer $l$ is uniquely determined
from the equality $c_yx=c_{yx}t^l$.

For each $x\in T/\gp t$, consider the set $X_x=X_1x\subseteq T/\gp t$
and the subproduct
$$
L_x=\left(\zvezda_{y\in X_x}G^{(c_y)}\right)*\gp z_\infty
$$
of the free product $L$. The group $L_x$ has a normal subgroup
$N_x=N_1^x\:=N_1^{\chi}$, where $\chi\in T$ is any representative of the
element $x\in T/\gp t$.

Let us show that the family of subproducts $\{L_x\ |\ x\in T/\gp t\}$
together with the subgroups $N_x\nin L_x$ satisfies the conditions
of Proposition 2. Indeed, conditions 1), 2), and 3) of Proposition 2
follow directly from the strong unique-product property of the group
$T/\gp t$. Condition $(**)$ for the pair $N_1\nin L_1$
follows from decomposition (3). Condition~$(**)$ for any other pair
$N_x\nin L_x$ also holds, because the groups $L_x$ and $L_1$ are
isomorphic and the isomorphism (the action of the element $x\in T$) maps
the subgroup $N_1$ onto the subgroup $N_x$ and each factor
$G^{(c_y)}$ of the group $L_1$ onto a subgroup $(G^{(c_{yx})})^{z^l}$ of
$L_x$.

Thus, the conditions of Proposition 2 hold.
Therefore, the natural mapping
$$
K_x=L_x/N_x\to K\:=L\Bigg/\!\gp{\!\!\!\gp{\bigcup_{y\in T/\gp t}N_y}\!\!\!}
$$
is injective.

The group $T$ acts on $K$ by automorphisms. Take the corresponding
semidirect product $T\semitimes K$ and consider its quotient by
the cyclic normal subgroup $\gp{zt^{-1}}$. The obtained group
$$
P=(T\semitimes K)/\gp{zt^{-1}}
$$
is the required quotient group of $\~G$.

Indeed, the group $G$ is embedded in $P$ as a subgroup:
$G=G^{(1)}\subseteq K \subseteq P$. According to the definition of
the action, we have
$
G^{(c_x)}=G^{c_x}.
$
Hence, the relations of the group $\~H_1$ (which are valid in $K$) and
the equality $t=z$ in $P$ give relation (2). Thus,
$P=\gp{T,G}$ is a quotient group of $\~G$ containing the subgroup $K_1$,
which, in its turn, contains any given countable group $S$. This completes
the proof of Theorem~2 in the case when the
group $\gp{\{t_i\}}$ is noncyclic.

\s 5.
Proof of Theorem 2. The splitting case

Now, let us prove Theorem 2 in the case when the group $\gp{\{t_i\}}$ is
cyclic. If presentation $(*)$ has the form $\~G=\pres<G*T|t=1>$, then the
group $\~G$ is the free product of two infinite groups $G$ and $T/\gp t$;
therefore, it is SQ-universal (see [LS77]). In what follows, we assume
that presentation~$(*)$ has a different form.

The unimodularity condition implies that
$\gp{\{t_i\}}=\gp t_\infty$, where $t=\prod t_i$.
Consider the group $R$ with the unimodular relative
presentation
$$
R=\pres<G,z|\prod_{i=1}^n g_iz^{k_i}=1>,
$$
where the exponents $k_i$ are determined from the equalities
$t_i=t^{k_i}$. We have $k_i\ne0$ and $\sum k_i=1$. It is known that
the element $z$ has infinite order in the group $R$ and
the equality $R=\gp z$ is valid
only if $n=1$ and $G=\gp{g_1}$ [CR01]. In this case,
Theorem 2 needs no proof.
In other cases, take an element~$r$ of the group $R$ not lying in $\gp z$,
but such that
$z^{kr}\ne z^l$ if the integers $k$ and $l$ are different.
It is easy to see that such an element exists. Indeed,
for all $x\in R$ let $k(x)$ be the non-negative integer defined by the
equality
$\gp z^x\cap\gp z =\gp{z^{k(x)}}$. There are three possible cases:
\item{1)}
$k(x)=0$ for some $x\in R$;
\item{2)}
$k(x)>1$ for some $x\in R$;
\item{3)}
$k(x)=1$ for all $x\in R$.

\enditem
In the first case we can take $r=x$.
In the second case we can take $r=z^{x^{-1}}$.
In the case 3), the cyclic subgroup $\gp z$ is normal in $R$
and does not coincides with its centraliser (because the index of this
centraliser in $R$ is at most two and any virtually cyclic torsion-free
group is cyclic);
in this case, any element of the centraliser of $z$ not lying in $\gp z$
can be taken as $r$.

Let
$\{t_{ij}\ ;\ i=1,2,\dots,\ j=1,\dots,1000\}$ be an infinite sequence of
elements of
the centraliser of $t$ in
$T$ such that all $t_{ij}^{\pm1}$ lie in different
cosets by the normal subgroup $\gp t$. Such a sequence exists, because
$T/\gp t$ is a nontrivial torsion-free group
and the index of the centraliser of $t$ in $T$ is at most two.
Take an arbitrary countable group $S=\{s_1,s_2,\dots\}$ and
put $K=R\times S$. Consider the group
$$
L=\left(K\zvezda_{z=t}T\right)
\Bigg/\!\!\gp{\!\!\!\gp{s_i\prod_{j=1}^{1000}rt_{ij}\ ;\
i=1,2,\dots}\!\!\!}.
$$
This presentation of the group $L$ satisfies the small cancellation condition
$C'(1/100)$ for free amalgamated products (see [LS77]). Therefore,
the natural mapping $S\to K\to L$ is injective. To complete the proof, it
remains to note that the group $L$ is a quotient of~$\~G=R\zvezda_{z=t}T$.


\s 6.
Proof of Theorem 3

\Lemma 2.
No infinite noncyclic group can be a union
of a finite number of its cyclic subgroups.

\Proof
Suppose that a group $G$ is a finite union of cyclic
subgroups. Obviously, such a group has the property that

\smallskip

{\sl any two infinite sets $X,Y\subseteq G$ contain a pair of
commuting elements $x\in X$, $y\in Y$}.

\smallskip

\noindent
It is well known that all infinite groups with this property
are abelian.%
\fn{%
This easily follows from a theorem of B.~Neumann [Neu76]
(the answer to a question of P.~Erd\H os):
\sl
The groups in which any infinite subset
contains a pair of different commuting elements are precisely
the groups with finite-index centres.
}
Since any subgroup and any quotient group of
the group $G$ is also a union of a finite number of cyclics,
it is sufficient to show that $G\not\iso\Z\oplus\Z_p$. But the group
$\Z\oplus\Z_p$ has infinitely many maximal cyclic subgroups:
$\gp{(1,1)}$, $\gp{(p,1)}$, $\gp{(p^2,1)}$,\dots; therefore, this
group cannot be a finite union of cyclics.

\Lemma 3.
Suppose that $l\ge2$ is an integer, $G_1,\dots,G_l$ are noncyclic
infinite groups, $u_1,\dots, u_s\in G_1*\dots*G_l$, and $S$ is
a countable group. Then the group $G_1*\dots*G_l$ has
a normal subgroup $N$ such that
\item{\rm1)}
the group $S$ embeds into $(G_1*\dots*G_l)/N$\rm;
\item{\rm2)}
$N\cap\gp{u_i,G_{i_1},\dots,G_{i_{l-1}}}=\1$ for all
$i\in\{1,\dots,s\}$ and $i_1,\dots,i_{l-1}\in\{1,\dots,l\}$.

\Proof
Let $U\subseteq G_1\cup\dots\cup G_l$ be the finite set consisting of
all coefficients of all words~$u_i$.
Lemma 2 implies that the set
$$
M_k=G_k\setminus\left(\bigcup_{u\in U}\gp u\right)
$$
is infinite for each $k\in\{1,\dots,l\}$.
Consider the countable set of words
$$
v_i=\prod_{j=1}^{\the\year}\prod_{k=1}^lg_{ijk} \in G_1*\dots*G_l,
$$
where $g_{ijk}\in M_k$ and all $g_{ijk}^{\pm1}$ are different.
Suppose that $S=\{s_1,s_2,\dots\}$ and put
$$
H=\left(G_1*\dots*G_l*S\right)/\nc{v_1s_1^{-1},v_2s_2^{-1},\dots}.
$$
Clearly, this presentation of the group $H$ satisfies the
small cancellation condition $C'(1/(100l))$ (see [LS77]). This
implies that the natural mapping $S\to H$ is injective, and each
word from the kernel $N$ of the natural epimorphism $G_1*\dots*G_l\to H$
contains coefficients from each of the sets~$M_k$. In particular,
$N\cap\gp{u_i,G_{i_1},\dots,G_{i_{l-1}}}=\1$, which completes the proof.

\medskip
Now, let us prove Theorem 3. Suppose that the word $w$ has the form
$u_1t^{\epsilon_1}\dots u_st^{\epsilon_s}$, where $u_i\in G_1*\dots*G_l$
and $\epsilon_i\in\{\pm1\}$. Note that, in the group $G_1*\dots*G_l$,
each nonidentity element~$u_i$ has infinite order and is transcendental
over each subgroup $G_{i_1}*\dots*G_{i_k}$ not containing this element
$u_i$.

Let us choose a normal subgroup $N\nin G_1*\dots*G_l$ according to Lemma 3
and put $G=(G_1*\dots*G_l)/N$.  Lemma~3 implies that the group $G$
contains any given countable group $S$ and the image of
each nonidentity element~$u_i$ in the group $G$ remains to be an element
of infinite order transcendental over each subgroups
$G_{i_1}*\dots*G_{i_k}$ not containing this element $u_i$. To complete the
proof, it remains to apply Proposition 1.


\s 7.
Proof of Proposition 1

We start with some simple facts.

\Lemma 4.
Let $u\in X*Y$ be an element of the free product of groups $X$ and
$Y$ and suppose that $Z$ is a subgroup of $Y$ such that the element $u$ is
algebraic (i.e., not transcendental) over $X*Z$. Then either
$u\in (X*Z)Y(X*Z)$
or $u$ has the form $x_1u'x_2$, where $x_1,x_2\in X*Z$ and the element
$u'\in X*Y$ has finite order.

We leave the proof of this elementary lemma to the readers as an easy
exercise.

\Lemma 5.
Suppose that $A$ is a nontrivial subgroup of a group $B$ and
$b\in B$. Then $b$ is transcendental over $A$ if and only if
$$
\gp{\{A^{b^i}\ ; \ i\in \Z\}}=\zvezda_{i\in \Z} A^{b^i}.
$$

\Proof
The ``only if" part is obvious. To prove the ``if" part, note that if
$u\in A*\gp{b}_\infty$ is a nontrivial relation between $A$ and $b$ in
the group $B$ and $a\in A\setminus\{1\}$, then $[a,u]$ is a nontrivial
relation between the groups $A^{b^i}$. Lemma 5 is proven.

Throughout this section, we assume that a subgroup $H$ of
a group $G$ satisfies conditions 1) and 2) from the definition of the strong
magnusianity, presentation (1) is unimodular, and all nonidentity
coefficients $g_i$ have infinite order in the group $G$. We have to
prove that the element $t$ is transcendental over $H$ in the group $\~G$.

\Lemma 6. If
$$
\gp{\{H^{t^i}\ ; \ i\in \Z\}}=\zvezda_{i\in \Z} H^{t^i}
\eqno{(4)}
$$
in the group $\~G$, then the element $t\in \~G$ is transcendental over~$H$.

\Proof
If the group $H$ is nontrivial, the assertion immediately follows from
Lemma 5. If $H=\{1\}$, then the lemma asserts only that the element $t\in
U$ has infinite order (if $w$ is not conjugate to $t^{\pm1}$); this fact
was proven in [CR01].

\smallskip

Put
$$
\=H=\zvezda_{i\in \Z} H_i\quad\hbox{and}\quad \=G=\=H\zvezda_{H_0=H}G,
$$
where $H_i$ are isomorphic copies of the group $H$.
Clearly, $\~G$ has the presentation
$$
\~G\iso\pres<\=G,t|w(t)=1,\ \{H_i^t=H_{i+1};\quad i\in\Z\}>.
$$
Let us move all letters $t^{\pm1}$ in the word $w$ to the left
through the coefficients
from $H$ by using the relations $H_it^{\pm1}=t^{\pm1}H_{i\pm1}$.
This reduces the presentation of $\~G$ to the form
$$
\~G\iso
\pres<\=G,t|\prod_{i=1}^p \=g_it^{k_i}=1,\
\{H_i^t=H_{i+1}\ ;\ i\in\Z\}>,
$$
where $k_i\in\Z\setminus\0$, $\sum k_i=1$, $\=g_i\in\=G$,
and each coefficient $\=g_i$ has the form
$$
\=g_i=\prod_{j=1}^s \=h_jf_j,
\quad\hbox{
where $s\ge1$, \
$f_j\in\{g_1,g_2,\dots\}\setminus H$, \
$\=h_j\in\=H\setminus H$ for $j\ne1$, \ and \ $\=h_1\in\=H$}.
$$
Here $f_j$, $\=h_j$, and $s$ depend on $i$.
Note that $f_1,\dots,f_s$ and
$\=h_2,\dots,\=h_s$ are transcendental over $H$ in $\=G$.

If $p=1$, then the first relation of $\~G$ can be
rewritten in the form
$$
t=u, \quad\hbox{where } u=\=g_1=\prod_{j=1}^s \=h_jf_j\in\=G,
$$
and the whole presentation can be rewritten in the form
$$
\~G\iso\pres<\=H \zvezda_{H_0=H} G|\{H_i^u=H_{i+1}\ ;\ i\in\Z\}>.
$$
By Lemma 6, it is sufficient to prove the injectivity of the natural
mapping $f\:\=H\to \~G$.

If $s>1$, i.e., $u\notin \=HG\=H$, then the homomorphism $f$ is
injective by virtue of Theorem 4.

Now, suppose that $s=1$, i.e., $u=\=h_1f_1$, where
$\=h_1\in \=H$, and the element $f_1\in G$ is
transcendental over $H$. In this case, the injectivity of the natural
mapping
$\=H$ to $\~G$
follows obviously from the decomposition
$$
\~G=G\zvezda_{K=L}(H*\gp t_\infty),
$$
where
$
G\supseteq K=\gp{f_1,H}=\gp{f_1}_\infty*H\iso
L=\gp{(\=h_1)^{-1}t,H}=\gp{(\=h_1)^{-1}t}_\infty * H
\subseteq H*\gp t_\infty
$.

We proceed to the case $p>1$.
Consider the following subgroups of $G*\gp{t}_\infty$: \
$G_i=t^{-i}Gt^i$,
$H_i=t^{-i}Ht^i$,
$$
\=H=\zvezda_{i=-\infty}^\infty H_i,\quad
K^{(m)}=\zvezda_{i=0}^m G_i, \quad\hbox{and}
\quad G^{(m)}=\=H\zvezda_{H_0*\dots*H_m}K^{(m)}.
\eqno{(5)}
$$
Consider all possible expressions of the relation $w=1$ in the form
$$
ct\prod_{i=1}^n b_it^{-1}a_it=1, \quad\hbox{where } a_i,b_i,c\in G^{(m)}.
\eqno{(6)}
$$
Among all such expressions we choose those in which $m$ is minimal; after
that, from all expressions with minimal $m$ we choose an expression with
minimal $n$.  For such a minimal expression (6), we have
\item{1)}
$n\ge 1$
(i.e., the length of this expression is strictly larger than one);
\item{2)}
$a_i\notin G^{(m-1)}$ and
$b_i\notin (G^{(m-1)})^t$;
\item{3)}
each coefficient $a_i$ is transcendental over $G^{(m-1)}$, and
each coefficient $b_i$ is transcendental over $(G^{(m-1)})^t$.

\noindent
The first property follows from the inequality $p>1$, which implies that
$m>0$ in any expression of length 1; therefore, $m$ can be decreased by
replacing all occurrences of elements of $G_m$ by fragments of the form
$t^{-1}gt$, where $g\in G_{m-1}$. The second property follows obviously
from the minimality of $n$ and $m$. To prove property~3),
we note that, in the normal form corresponding to the decomposition
$$
G^{(m)}=X*Y,\quad\hbox{
where $X=\left(\zvezda_{j\ne m}H_j\right)*G_0*\dots G_{m-1}$
and $Y=G_m$,
}
$$
each $Y$-syllable of each coefficient $a_i$ lies
in $(H\{f_1,\dots,f_p\}H)^{t^m}$, and apply Lemma 4 putting $Z=H_m$.
The transcendence of $b_i$ over $(G^{(m-1)})^t$ is
established similarly.

Now, suppose that the symbols $H_i$ and $G_i$ denote abstract isomorphic
copies of the groups $H$ and $G$ and the groups $\=H$, $K^{(m)}$, and
$G^{(m)}$ are defined by formulae (5). Consider the following
presentation of the group $\~G$.
$$
\~G= \gp{G^{(m)},t\ \biggm|\
ct\prod_{i=1}^n b_ia_i^t=1,\
\left\{G_i^t=G_{i+1}\ ;\quad i\in\{0,\dots m-1\}\right\},\
\{H_i^t=H_{i+1}\ ;\quad i\in\Z\}\
}.
\eqno{(7)}
$$
To complete the proof of Proposition 1, it remains to note that
properties 1) and 3) of presentation (7) imply the injectivity
of the natural mapping $G^{(m)}\to\~G$ by virtue of the following theorem.

\Th \rm{([Kl93], see also [Fer96])}.
Suppose that $M$ and $N$ are
isomorphic subgroups of a group $L$, $\phi:M\to N$ is an isomorphism,
$n\ge 1$, $a_1,\dots, a_n$ are elements of $L$
transcendental over $M$,
$b_1,\dots, b_n$ are
elements of $L$ transcendental over $N$, and $c\in L$. Then the system
of equations
$$
\left\{\eqalign{
&x^{-1}gx=g^\phi,\quad g\in M, \cr
&cx\prod_{i=1}^n b_ix^{-1}a_ix=1
}\right.
\eqno{({**}*)}
$$
is solvable over $L$, i.e., the natural mapping
$L\to\pres<L,x|({**}*)>$ is injective.

Applying this theorem to $L=G^{(m)}$ and $M=G^{(m-1)}$,
we see that the natural mapping $\=H\subset G^{(m)}\to\~G$ is
injective and the element $t\in\~G$ is transcendental over $H$ (by
Lemma 6).


\s 8.
Proof of Theorem 4

Let us reformulate Theorem 4 in terms of equations over
groups. Recall that an {\it equation over a group $G$ with unknown
\({\rm or} variable\) $t$} is a formal expression of the form
$$
g_1t^{\epsilon_1}g_2t^{\epsilon_2}\dots g_nt^{\epsilon_n}=1,
\eqno{(8)}
$$
where $g_i\in G$ and $\epsilon_i\in\Z$.
Equation $(8)$ is said to be
{\it solvable over the group $G$} if there exists an overgroup
$\~G$ of $G$ and an element
$\~t\in\~G$ (called a solution to equation $(*)$) such
that
$g_1\~t^{\epsilon_1}g_2\~t^{\epsilon_2}\dots g_n\~t^{\epsilon_n}=1$ in
$\~G$. The notion of {\it solvability of a
system of equations with several unknowns over a group $G$} is
defined similarly.

\Lemma 7.
Suppose that $G=A\zvezda_C B$ is the amalgamated free product of
groups $A$ and $B$ with amalgamated subgroup $C$,
$v=b_0a_0\dots b_ma_mb_{m+1}\in G$, $\phi$ is an automorphism of $B$,
and
$$
\^G=\gp{A\zvezda_C B\ \biggm|\ \{b^v=b^\phi\ |\ b\in B\}}.
$$
Then the injectivity of the natural mappings $A\to \^G\leftarrow B$ is
equivalent to the solvability of the following system of equations with
unknowns $t$ and $x$ over the group $G$:
$$
\cases{
b^{-x}b^\phi=1 \quad
\hbox{ for }
b\in B\setminus\1,\cr
[t,c]=1 \quad\ \
\hbox{ for }
c\in C\setminus\1,\cr
x^{-1}b_0a_0^t\dots b_ma_m^tb_{m+1}=1.
}
\eqno{(9)}
$$

\Proof
If $\~t,\~x\in\~G\supseteq G$ is a solution to system (9), then
the mapping $a\mapsto a^{\~t}$, $b\mapsto b$ extends to a homomorphism
$\^G \to \~G$ injective on
$A$ and on $B$.
Conversely, suppose that the natural mappings
$
A\to \^G\leftarrow B
$
are injective and consider an isomorphic copy
$$
\=G=\gp{\=A\zvezda_{\=C}\=B\ \biggm|\ \{\=b^{\=v}=\=b^{\=\phi}\ |\
b\in B\}}
$$
of the group $\^G$.
It is easy to see that
the elements $\~t$ and
$\~x=b_0a_0^{\~t}\dots b_ma_m^{\~t}b_{m+1}=\=v$ \ of the HNN-extension
$$
\pres<G\zvezda_{B=\=B}\=G, \~t | \{a^{\~t}=\=a\ ;\ a\in A\}>
$$
are solutions to system (9). This proves Lemma 7.

\smallskip

This lemma shows that Theorem 4 is implied by the following
theorem.

\Th 4$'$.
Suppose that $G$ is a group with subgroups $C$, $B$, and $B^\phi$,
$\phi\:B\to B^\phi$ is an isomorphism, $m\ge1$, $b_i,a_i\in G$, and the
elements $a_0,\dots, a_m$ and $b_1,\dots, b_m$ are transcendental over
$C$.  Then system (9) is solvable over $G$.

To prove Theorem 4$'$, we need Howie's lemma.  For simplicity, we
formulate this lemma in a special case related to system (9).

Consider a map (tessellation) on an oriented two-dimensional sphere.
The corners of this map are labelled by elements of the group $G$. The
edges are directed (the directions are shown by arrows on the figures) and
labelled by the variables $t$ and $x$.

The {\it label of a vertex} in such a situation is defined as the product
of the labels of all corners near this vertex listed clockwise. The label
of a vertex is an element of $G$ defined up to conjugation. For example,
the label of the vertex shown in Figure 2 is $a_0^{-1}ca_0^2c'$.

\goodbreak
\bigskip
\centerline{\input 2.PIC}
\nobreak%
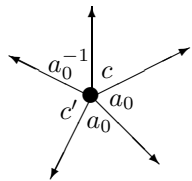
\centerline{Fig. 2}%
\goodbreak
\bigskip

To obtain the {\it label of a face}, we should walk along its boundary
anticlockwise and write down the labels of all its corners and edges;
the label of an edge should be written as its inverse if we
walk through it against the arrow.
The label of a face is an element of the group $G*F(t,x)$ (the free
product of $G$ and the free group with basis $\{t,x\}$),
defined up to a cyclic permutation.
For example, the label of the upper left face shown in Figure 3 is
$b^{-x}b^\phi$.

Such a labelled map is called a {\it spherical Howie diagram} (or simply
{\it diagram}) over system (9) if
\item{1)} one vertex is
  distinguished and called the {\it exterior vertex},
  the other vertices are called {\it interior};
\item{2)}
  the label of each interior vertex is the identity element of $G$;
\item{3)}
  the label of each face is either the left-hand side of an equation of
  system (9) or the word inverse to the left-hand side of an equation
  of system (9); all possible types of faces are shown in Figure 3, where
  the letters $b$ and $c$ denote any nonidentity elements of the groups
  $B$ and $C$, respectively, and the edges not labelled by $x$ are assumed to
  be labelled by $t$ (the numbers written outside the faces should be ignored
  for a while).

\goodbreak
\bigskip
\centerline{\input 3.PIC}
\nobreak%
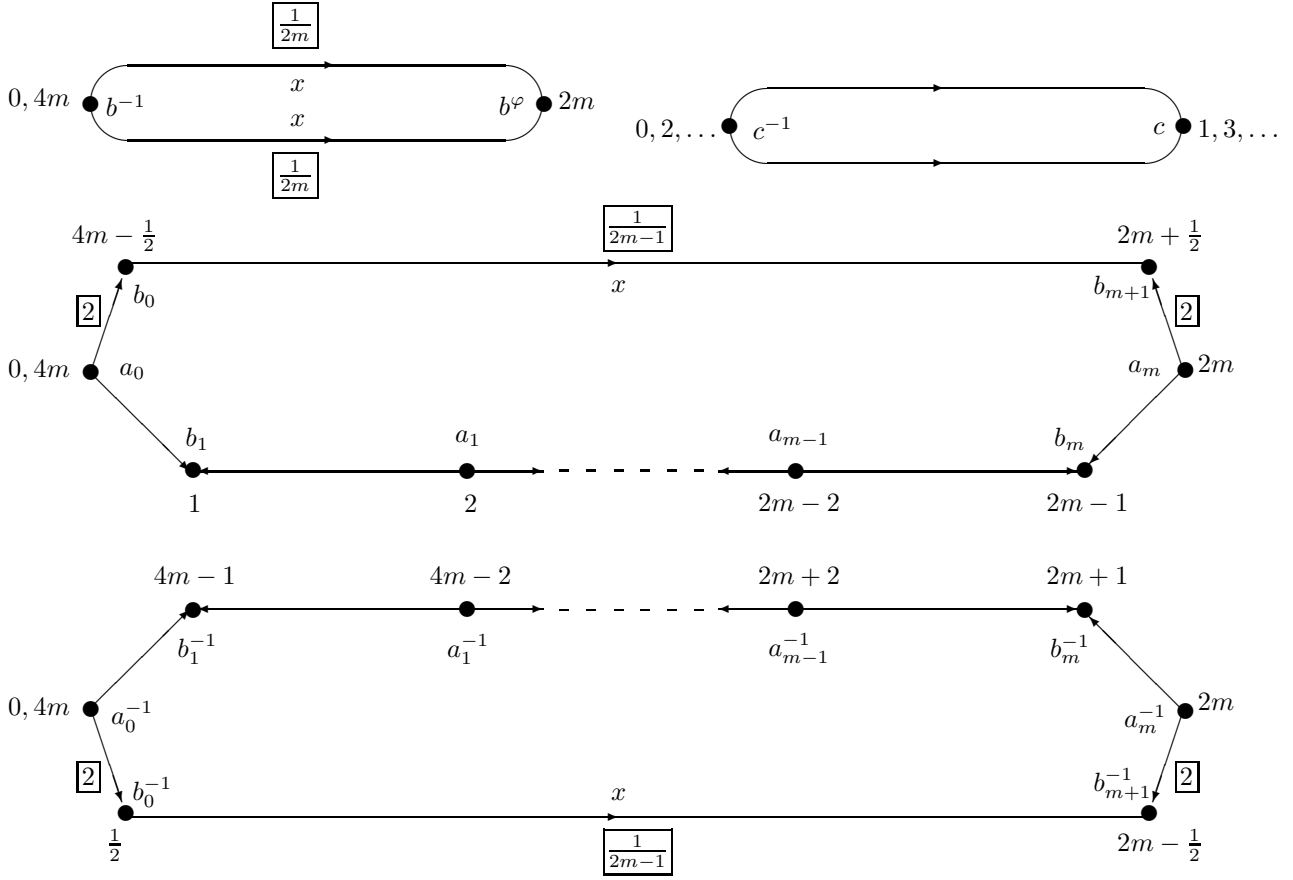
\centerline{Fig. 3}%
\goodbreak
\bigskip

A Howie diagram is called {\it reduced} if it has no
edge $e$ such that the two faces containing $e$ are
different and their labels written starting with the edge $e$ are
mutually inverse; such a pair of faces with a common edge is called a {\it
reducible pair}.

\Lemma 8 {\rm[How83]}.
System (9) is not solvable over the group $G$
if and only if there exists a spherical
diagram over this system such that the label of its exterior vertex is
not the identity element of $G$. Any such diagram with minimal number
of faces is reduced.

We say that a diagram over system (9) is {\it strongly reduced} if it is
reduced and different faces labelled by $b^{-x}b^\phi$ or $[c,t]$ have
no common edges.

\Lemma 9.
Any spherical diagram with minimal number of faces such that the label of
its exterior vertex is a nonidentity element is strongly reduced.

\Proof
Indeed, if some diagram has a pair of
faces labelled by, e.g., $b^{-x}b^\phi$ and $(b')^{-x}(b')^\phi$ and
having a common edge, then either this pair is a reducible pair or we can
erase the common edge multiplying the labels of the corresponding corners
(Fig. 4) and obtain a diagram fewer faces and the same
label of the exterior vertex. This implies the nonminimality of the
initial diagram and proves the lemma.

\goodbreak
\bigskip
\centerline{\input 4.PIC}
\nobreak%
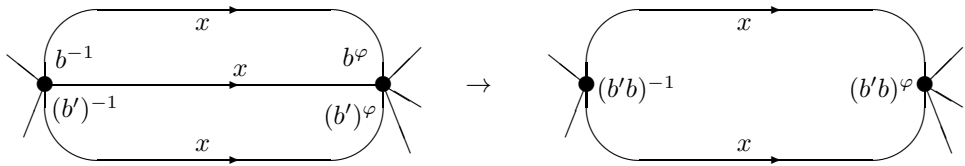
\centerline{Fig. 4}%
\goodbreak
\bigskip

Now, suppose that, on the contour of each face
$D$ of a map on the sphere, there is a moving point (a {\it car})
$\alpha_D$. The car $\alpha_D$ moves continuously anticlockwise (i.e.,
the interior of the face $D$ remains on the left from the car)
without U-turns, stops, and ``infinite decelerations", that is, it
covers each edge in a finite time.
We call such a motion {\it regular}.

If the number of cars being at a moment $\tau\in\R$ at
a point $p$ of the sphere equals the multiplicity of this
point (in other words, either two cars simultaneously pass the same
internal point $p$ of an edge at the moment
$\tau$ or
there are $k$ cars at a vertex $p$ of degree $k$ at the moment $\tau$),
then we say that a {\it complete collision} occurs at the point $p$ at the
moment $\tau$; the point $p$ is called a {\it point of complete
collision}. Figure 5 shows complete collisions on an edge and at a vertex
of degree three.

\goodbreak
\bigskip
\centerline{\input 5.PIC}
\nobreak%
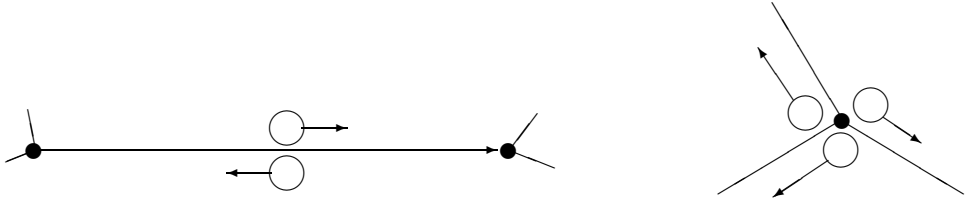
\centerline{Fig. 5}%
\goodbreak
\bigskip

\Lemma 10
{\rm[Kl93] (see also [FeR96])}.
Any regular motion on a sphere has at least two points of complete
collision.

Now, take a spherical diagram over presentation
(9) as a map and consider the following regular motion on this map:

\item{a)}
  the car going around a face labelled by $b^{-x}b^\phi$ moves
  anticlockwise uniformly with speed $1\over{2m}$ edges per unit time and
  visits the corner labelled by $b^{-1}$ at time zero.

\item{b)}
  the car going around a face labelled by $[t,c]$ moves anticlockwise
  uniformly with unit speed (one edge per unit time) and visits the
  corner with label $c^{-1}$ at time zero;

\item{c)}
  the car going (anticlockwise) around a face labelled by
  $x^{-1}b_0a_0^t\dots b_ma_m^tb_{m+1}$
  is at the corner labelled by $a_0$ at time zero; then it goes along
  $2m$ edges labelled by $t$ with unit speed and arrives at the moment $2m$
  in the corner labelled by $a_m$; then, the car goes along the next edge
  labelled by $t$ with speed 2 and the edge labelled by
  $x$ with speed $1\over{2m-1}$; after that, it moves along the edge
  labelled by $t$ with speed 2 and returns at the moment $4m$ to the
  initial corner labelled by $a_0$; then, the motion is repeated with
  period~$4m$;

\item{d)}
  the car moving (anticlockwise) around a face labelled by
  $(x^{-1}b_0a_0^t\dots b_ma_m^tb_{m+1})^{-1}$ is at the corner labelled
  by~$a_0^{-1}$ at time zero; then it goes along the first edge labelled
  by $t$ with speed 2 and along the next edge labelled by~$x$ with
  speed $1\over{2m-1}$; after that, the car goes along the edge labelled
  by $t$ with speed 2 and arrives at the moment~$2m$ in the corner
  labelled by $a_m^{-1}$; then, the car goes along the next $2m$ edges
  labelled by $t$ with unit speed and returns at the moment $4m$ to the
  initial corner labelled by $a_0^{-1}$; then the motion is repeated with
  period $4m$.

\noindent
This motion is regular and periodic with period $4m$
(on the faces labelled by $[t,c]$, the minimal period is 2).
Figure~3 shows the detailed schedule of this motion during the interval
$0\le\tau\le 4m$; the boxed numbers near edges indicate the speed of the
car on these edges (by default, the speed is 1).

\Lemma 11.
For the motion on a strongly reduced diagram over system
(9) described above, the complete collisions can occur only
at the exterior vertex.

\Proof

\noindent{\bb A collision on an edge}
labelled by $t$ at time $\tau$ means that, at this moment, the direction
of the motion of one of the cars coincides with the direction the edge,
while the direction of the motion of the other colliding car is opposite
to the direction of this edge (Fig. 5, left). But the schedule of the
motion is such that, at each moment $\tau$, either all cars being on
the edges labelled by $t$ move in the direction of the edge (this is so
when the integer part of $\tau$ is even) or all cars being on edges
labelled by $t$ move in the direction opposite to the direction of the
edge (this is so when the integer part of $\tau$ is odd).

A similar argument shows that cars cannot collide on an edge labelled
by $x$: during the time intervals $[0, 2m]+4m\Z$, all
cars being on edges labelled by $x$ move in the direction of the edge,
and during the time intervals $[2m,4m]+4m\Z$,
all cars being on edges labelled by $x$ move in the direction
opposite to the direction of the edge.

Thus, collisions can occur only at vertices.

\noindent{\bb A complete collision at the start or the end vertex of an edge
labelled by $x$}
cannot occur, because the diagram is strongly reduced and, hence, each
edge labelled by $x$ separates a face labelled by $b^{-x}b^\phi$ and a
face labelled $(x^{-1}b_0a_0^t\dots b_ma_m^tb_{m+1})^{\pm1}$ (see Fig. 6).
This means that one car visits the start vertex of an edge labelled by
$x$ at the moments $4m\Z$, and another car visits this vertex at the
moments $4m\Z\pm{1\over2}$. Thus, a complete collision cannot occur. For
similar reasons, a complete collision cannot occur at the end vertex of an
edge labelled by $x$: one car visits such a vertex at the moments
$4m\Z+2m$ and another, at the moments $4m\Z+2m\pm{1\over2}$.

\goodbreak
\bigskip
\centerline{\input 6.PIC}
\nobreak%
\centerline{Fig. 6}%
\goodbreak
\bigskip

\noindent{\bb A complete collision at an interior vertex being
neither the start nor the end of an edge labelled by $x$} cannot
occur too. Such vertices are visited only at integer moments of time
and, at every such moment $\tau$, each car being at such a
vertex is either at a corner labelled by $c\in C\setminus\1$ or at a
corner labelled by
$d\in\{a_0^{\pm1},b_1^{\pm1},\dots,b_m^{\pm1},a_m^{\pm1}\}$.  Moreover, if
at the moment $\tau$ a car is at a corner labelled by $d$ and,
simultaneously, another car is at a corner labelled by $d'$, where
$d,d'\in\{a_0^{\pm1},b_1^{\pm1},\dots,b_m^{\pm1},a_m^{\pm1}\}$, then
either $d'=d$ or $d'=d^{-1}$.%
\fn{%
The latter can occur only at the moments $\tau\in\{0,2m\}+4m\Z$. Thus,
$d=a_0^{\pm1}$ or $d=a_m^{\pm1}$. This makes it possible to slightly weaken
the conditions of Theorems 4 and 4$'$, but we do not need such a weakening
here.
}
This means that a vertex of complete collision occurring at the
moment $\tau$ has label of the form $\prod z_i$, where $z_i\in
\{d,d^{-1}\}\cup C\setminus\1$.  Since the diagram is reduced, $d$ and
$d^{-1}$ cannot be neighbours in the sequence $(z_i)$. Since the diagram
is strongly reduced, two elements from $C\setminus\1$ cannot be neighbours
in this sequence either. Therefore, the label $\prod z_i$ of the vertex of
complete collision cannot be the identity element of $G$, because the
element $d$ is transcendental over $C$. Hence, the vertex of complete
collision cannot be an interior vertex of the diagram. Figure 2 shows a
hypothetical vertex at which a complete collision occurs at the moment
$\tau=0$. This vertex cannot be interior, because $a_0^{-1}ca_0^2c'\ne1$
in the group $G$.  This completes the proof of Lemma 11.

Theorem 4$'$ follows easily from lemmata 8--11. Indeed, suppose
that system (9) is unsolvable. Then, by Lemma~8, there exists a diagram
over this system. By Lemma 9, this diagram can be assumed strongly
reduced. By Lemma~11, it admits a regular motion with at most
one point of complete collision, which contradicts Lemma 10.
This proves Theorems 4$'$ and 4.


\s{\rm REFERENCES}

\item{[B84]}
Brodskii S.D.
Equations over groups and one-relator groups
{// Sib. Mat. Zh.} 1984. {T.25}. no.2. P.84--103.

\item{[Kl94]}
Klyachko Ant.A.
The Kervaire--Laudenbach Conjecture and Equations over Groups,
{Cand. Sci. Dissertation},
{Moscow: MSU,} 1994.

\item{[Kl05]}
Klyachko Ant.A.
The Kervaire--Laudenbach conjecture and presentations of simple groups
{// Algebra i Logika}. 2005. {T. 44}. {no.4}. P. 399--437.
See also
{arXiv:math.GR/0409146}.

\item{[Kl06a]}
Klyachko Ant.A.
How to generalize known results on equations over groups
{// Mat. Zametki}. 2006. {T.79}. no.3. P.409--419.
See also
{arXiv:math.GR/0406382}.

\item{[Kl06b]}
Klyachko Ant.A.
Free subgroups of one-relator relative presentations
{// Algebra i Logika} (to appear).
See also
{arXiv:math.GR/0510582}.

\item{[LS77]}
Lyndon R.C., Schupp P.E.
{Combinatorial Group Theory},
Springer-Verlag, Berlin/Heidelberg/New~York, 1977.

\item{[Lo86]}
Lossov K.I.
SQ-universality of free products with finite amalgamated subgroups
{// Sib. Mat. Zh.} 1986. {T.27}. no.6. P.128--139.

\item{[Ols95]}
Olshanskii A.Yu.
SQ-universality of hyperbolic groups
{// Mat. Sbornik}. 1995. {T.186}. no.8. P.119--132.

\item{[AMO06]}
Arzhantseva G., Minasyan A., Osin D.
The SQ-universality and residual properties of relatively hyperbolic groups
{// arXiv:math.GR/0601590}.

\item{[BaPr78]}
Baumslag B., Pride S.
Groups with two more generators than relators
{// J. London Math. Soc.} 1978. {V.17}. P.425--426.

\item{[Bu05]}
Button J.O.
Large mapping tori of free group endomorphisms
{// arXiv:math.GR/0511715}.

\item{[CR01]}
Cohen M.M., Rourke C.
The surjectivity problem for one-generator, one-relator extensions of
torsion-free groups
{// Geometry \& Topology}. 2001. {V.5}. P.127--142.
See also
{arXiv:math.GR/0009101}.

\item{[Ed84]}
Edjvet M.
Groups with balanced presentations
{// Arch. Math.} 1984. {V.42}. no.4. P.311--313.

\item{[FeR96]}
Fenn R., Rourke C.
Klyachko's methods and the solution of equations over torsion-free groups
{// L'Enseignment Math\'ematique.} 1996. {T.42}. P.49--74.

\item{[FoR05]}
Forester M., Rourke C.
Diagrams and the second homotopy group
{// Comm. Anal. Geom.} 2005. {V.13}. P.801-820.
see also
arXiv:math.AT/0306088.

\item{[Gr83]}
Gromov M.
Volume and bounded cohomology
{// Inst. Hautes Etudes Sci. Publ. Math.}
1982. No. 56. P.5--99. (1983).

\item{[How83]}
Howie J.
The solution of length three equations over groups
{// Proc. Edinburgh Math. Soc.} 1983. {V.26}. P.89--96.

\item{[How98]}
Howie J.
Free subgroups in groups of small deficiency
{// J. Group Theory}. 1998. V.1. no. 1.  P.95--112.

\item{[Kl93]}
Klyachko Ant.A.
A funny property of a sphere and equations over groups
{// Comm. Algebra}. 1993. {V.21}. P.2555--2575.

\item{[La05]}
Lackenby M.
A characterisation of large finitely presented groups
{// J. Algebra}. 2005. {V.287}. {no.2}. P. 458--473.
See also
{arXiv:math.GR/0403129}.

\item{[Neu76]}
Neumann B.H.
A problem of Paul Erd\H os on groups
{// J. Austral. Math. Soc. Ser. A.} 1976. {V.21}. {no.4}. P.467--472.

\item{[Neu73]}
Neumann P.M.
The SQ-universality of some finitely presented groups
{// J. Austral. Math. Soc.} 1973. {V.16}. P.1--6.

\item{[OlOs06]}
Olshanskii A.Yu., Osin D.V.
Adding high-powered relations to large groups:
A short proof of Lackenby's result
{// arXiv:math.GR/0601589}.

\item{[P88]}
Promyslow S.D.
A simple example of a torsion free nonunique product group
{// Bull. London Math. Soc.} 1988. {V.20}. P.302--304.

\item{[RS87]}
Rips E., Segev Y.
Torsion free groups without unique product property
{// J. Algebra} 1987. {V.108}. P.116--126.

\item{[SaSc74]}
Sacerdote G.S., Schupp P.E.
SQ-universality in HNN groups and one relator groups
{// J. London Math. Soc.} 1974. {V.7}. P.733--740.

\item{[St\"o83]}
St\"ohr R.
Groups with one more generator than relators
{// Math. Z.} 1983. {V.182}.  no. 1. P.45--47.

\end